\documentclass[11pt]{article}
\usepackage{latexsym,amsfonts,amssymb,amsmath,amsthm}

\usepackage{color,comment}
\usepackage{mathtools}

\begin{document}

\newcommand{ \bl}{\color{blue}}
\newcommand {\rd}{\color{red}}
\newcommand{ \bk}{\color{black}}
\newcommand{ \gr}{\color{OliveGreen}}
\newcommand{ \mg}{\color{RedViolet}}

\newcommand{\norm}[1]{||#1||}
\newcommand{\normo}[1]{|#1|}
\newcommand{\secn}[1]{\addtocounter{section}{1}\par\medskip\noindent
 {\large \bf \thesection. #1}\par\medskip\setcounter{equation}{0}}

\newcommand{\secs}[1]{\addtocounter {section}{1}\par\medskip\noindent
 {\large \bf  #1}\par\medskip\setcounter{equation}{0}}
\renewcommand{\theequation}{\thesection.\arabic{equation}}

\setlength{\baselineskip}{16pt}

\newtheorem{theorem}{Theorem}[section]
\newtheorem{lemma}{Lemma}[section]
\newtheorem{proposition}{Proposition}[section]
\newtheorem{definition}{Definition}[section]
\newtheorem{example}{Example}[section]
\newtheorem{corollary}{Corollary}[section]

\newtheorem{remark}{Remark}[section]

\numberwithin{equation}{section}

\def\p{\partial}
\def\I{\textit}
\def\R{\mathbb R}
\def\C{\mathbb C}
\def\u{\underline}
\def\l{\lambda}
\def\a{\alpha}
\def\O{\Omega}
\def\e{\epsilon}
\def\ls{\lambda^*}
\def\D{\displaystyle}
\def\wyx{ \frac{w(y,t)}{w(x,t)}}
\def\imp{\Rightarrow}
\def\tE{\tilde E}
\def\tX{\tilde X}
\def\tH{\tilde H}
\def\tu{\tilde u}
\def\d{\mathcal D}
\def\aa{\mathcal A}
\def\DH{\mathcal D(\tH)}
\def\bE{\bar E}
\def\bH{\bar H}
\def\M{\mathcal M}
\renewcommand{\labelenumi}{(\arabic{enumi})}

\def\disp{\displaystyle}
\def\undertex#1{$\underline{\hbox{#1}}$}
\def\card{\mathop{\hbox{card}}}
\def\sgn{\mathop{\hbox{sgn}}}
\def\exp{\mathop{\hbox{exp}}}
\def\OFP{(\Omega,{\cal F},\PP)}
\newcommand\JM{Mierczy\'nski}
\newcommand\RR{\ensuremath{\mathbb{R}}}
\newcommand\CC{\ensuremath{\mathbb{C}}}
\newcommand\QQ{\ensuremath{\mathbb{Q}}}
\newcommand\ZZ{\ensuremath{\mathbb{Z}}}
\newcommand\NN{\ensuremath{\mathbb{N}}}
\newcommand\PP{\ensuremath{\mathbb{P}}}
\newcommand\abs[1]{\ensuremath{\lvert#1\rvert}}

\newcommand\normf[1]{\ensuremath{\lVert#1\rVert_{f}}}
\newcommand\normfRb[1]{\ensuremath{\lVert#1\rVert_{f,R_b}}}
\newcommand\normfRbone[1]{\ensuremath{\lVert#1\rVert_{f, R_{b_1}}}}
\newcommand\normfRbtwo[1]{\ensuremath{\lVert#1\rVert_{f,R_{b_2}}}}
\newcommand\normtwo[1]{\ensuremath{\lVert#1\rVert_{2}}}
\newcommand\norminfty[1]{\ensuremath{\lVert#1\rVert_{\infty}}}
\newcommand{\ds}{\displaystyle}

\title{Uniqueness and nonlinear stability of positive entire solutions in parabolic-parabolic chemotaxis models with logistic source on bounded heterogeneous environments}

\author{ Tahir Bachar Issa
\thanks{tahirbachar.issa@sjsu.edu, Department of Mathematics and Statistics,\\
San Jose State University, San Jose, CA, 95192, USA.}
}

\date{}
\maketitle

\noindent {\bf Abstract.} {This paper studies the asymptotic behavior of solutions of the parabolic-parabolic chemotaxis model with logistic-type sources in heterogeneous bounded domains:
\begin{equation*}
 \label{u-v-eq00}
\begin{cases}
u_t=\Delta u-\chi\nabla\cdot (u \nabla v)+u\Big(a_0(t,x)-a_1(t,x)u-a_2(t,x)\int_{\Omega}u\Big),\quad x\in \Omega\cr
\tau v_t=\Delta v-\lambda v +\mu u,\quad x\in \Omega \cr
\frac{\p u}{\p \nu}=\frac{\p v}{\p \nu}=0,\quad x\in\p\Omega.
\end{cases}\qquad(\ast)
 \end{equation*}
\noindent
We find parameter regions in which the system has a unique positive entire solution, which is globally asymptotically stable.
{More precisely under suitable assumptions on the model's parameters, the system has a unique {positive entire solution} $(u^*(t,x),v^*(t,x))$ such that for any 
$u_0 \in C^0(\bar{\Omega}),$ $v_0 \in W^{1,\infty}(\bar{\Omega})$ with $u_0,v_0\ge 0$ and
$u_0\not\equiv 0$,  the  global classical solution
$(u(t,x;t_0,u_0,v_0)$,  $v(t,x;t_0,u_0,v_0))$ of $(\ast)$ satisfies
$$
  \lim_{t \to \infty}\Big(\sup_{t_0 \in \mathbb{R}}\|u(t+t_0,\cdot;t_0,u_0,v_0)-u^*(t+t_0,\cdot)\|_{C^0(\bar\Omega)}+\sup_{t_0 \in \mathbb{R}}\|v(t+t_0,\cdot;t_0,u_0,v_0)-v^*(t+t_0,\cdot)\|_{C^0(\bar\Omega)}\Big)=0.
$$
}
\medskip

\noindent {\bf Keywords.} {Fully {parabolic chemotaxis model}, entire solutions, uniqueness and stability, comparison principle.}

\medskip

\noindent {\bf 2010 Mathematics Subject Classification.} 35A09, 35B08, 35B09, 35B40

\section{Introduction and the statements of the main results}
\label{S:intro}
Chemotaxis is the process of directed movement of organisms towards a higher or lower concentration of a specific chemical substance. It plays a crucial role in various biological phenomena, such as immune system response, embryo development, tumor growth, population dynamics, and gravitational collapse \cite{ISM04, DAL1991}. In 1970, Keller and Segel proposed a celebrated mathematical model to describe the chemotaxis phenomenon \cite{KS1970, KS71}. {Since these interesting works of Keller and Segel}, many variants of chemotaxis models have been proposed and their dynamics have been studied intensively mainly on homogeneous environments \cite{ASV2009, JDITBRB, HVJ1997a,HMVJ1996,HeSh,HW2001,ITBRS17,ITBWS16,ITBRBSWS2021,JaW92,HIKWS2023b,HIKWS2023a,kuto_PHYSD, PaHi, TW07,Zhe} and \cite{LiKeWa,LiMu,NT1995,NT2001,NT13,NTa15,RBSWS17a,T04,TJiCMAJ,TW07,Win2010,W2011b,Win2014}. In particular, we refer to the survey papers \cite{NBYTMW05, THKJP09,  H03}. In reality, the underlying environments of many biological systems are subject to various spatial and temporal variations. Studying chemotaxis models in heterogeneous environments is of both biological and mathematical importance. Recently, chemotaxis models in heterogeneous environments have been introduced, and their dynamics have been investigated in various studies \cite{ITBWS16, ITBWS17a, ITBWS17b, ITBWS2020,RBSWS2018,RBSWS2018b}. However, important dynamical questions, including the uniqueness and global stability of positive entire solutions, remain open for parabolic-parabolic chemotaxis models in heterogeneous environments.

The objective of this paper is to investigate the uniqueness and stability of positive entire solutions of the full parabolic-chemotaxis model with time and {space dependent logistic source} 
\begin{equation}
 \label{u-v-eq00}
\begin{cases}
u_t=\Delta u-\chi\nabla\cdot (u \nabla v)+\big(a_0(t,x)-a_1(t,x)u-a_2(t,x)\int_{\Omega}u\big)u,& x\in \Omega,\cr
\tau v_t=\Delta v-\lambda v +\mu u,& x\in \Omega, \cr
\frac{\p u}{\p \nu}=\frac{\p v}{\p \nu}=0,& x\in\p\Omega,
\end{cases}
 \end{equation}
where $\Omega \subset \mathbb{R}^n(n\geq 1)$ is a bounded domain with smooth boundary $\p\Omega$, and $\nu$ is the outward normal unit vector to $\p\Omega$. $u(x,t)$ represents the local population density
of a  mobile species, $v(x,t)$ is the local population density of the chemical substance created by the mobile species, $\chi \in \mathbb{R}$ represents
the chemotactic sensitivity, $\tau$ is a positive constant related to the diffusion rate of the chemical substance, { $\lambda$ is a positive constant that represents the degradation rate of the chemical substance, $\mu$ is a positive constant that presents the rate
at which the mobile species produces the chemical substance.
{The first equation in \eqref{u-v-eq00} contains the term $\big(a_0(t,x)-a_1(t,x)u-a_2(t,x)\int_\Omega u\big)u$, which describes how individuals of the species compete for resources while also cooperate to survive. The positive function $a_0(t,x)$ leads to exponential growth when the population density is low, while the term $a_1(t,x)u$ represents local competition. 
The nonlocal term $a_2(t,x)\int_\Omega u$ accounts for the effects of the species' total mass on its growth. If $a_2(t,x)>0$, the nonlocal term limits the growth of the mobile species while it leads to the population growth when $a_2(t,x)<0$. In the latter case, the mobile species compete locally but cooperate globally.  In 
 system \eqref{u-v-eq00}, we impose the homogeneous Neumann boundary conditions to indicate that there is no flux across the boundary of the habitat.} 
 }
 
 The paper \cite{ITBWS2020} has established the {global} existence and boundedness of classical solutions, pointwise persistence, and the existence of positive entire solutions for system \eqref{u-v-eq00}. Our aim in the current paper is to find an appropriate parameters range for which the positive entire solution of system \eqref{u-v-eq00} established in \cite{ITBWS2020} is unique and globally stable with respect to positive perturbations.

\subsection{Known results on system \eqref{u-v-eq00}}

{In the following, we recall some results established in \cite{ITBWS2020} on the global boundedness and persistence of classical solutions and the existence of at least one positive entire solution of system \eqref{u-v-eq00}}. To this end, we first introduce some notations and definitions.  Throughout the paper, we  put
\vspace{-0.05in} \begin{equation*}
 \label{a-i-sup-inf-eq1}
a_{\inf}=\inf _{ t \in\RR,x \in\bar{\Omega}}a(t,x),\ a_{\sup}=\sup _{ t \in\RR,x \in\bar{\Omega}}a(t,x), \ a_{\inf}(t)=\inf _{x \in\bar{\Omega}}a(t,x), \ \text{and}\ a_{\sup}(t)=\sup _{x \in\bar{\Omega}}a(t,x)
 \vspace{-0.05in}\end{equation*}
{for any $a\in L^{\infty}(\mathbb{R}\times\bar{\Omega})$ and $t\in\mathbb{R}$. For a given $z\in\mathbb{R}$,  we use standard notations $z_+=\max\{0,z\}$ and $z_-=\max\{0,-z\}$.}  We recall  definition of the classical solution of system \eqref{u-v-eq00}.
\begin{definition}
\label{global existence}
Given $t_0\in\RR$, $u_0\in C(\bar\Omega)$,  and $v_0\in W^{1,\infty}(\Omega)$ with $u_0> 0$ and  $v_0\geq 0$, a classical solution
$(u(t,x),v(t,x))$  is denoted by $(u(t,x;t_0,u_0,v_0),v(t,x;t_0,u_0,v_0))$ if it is defined on $[t_0,t_0+T)$ for some $T\in(0,\infty]$
and satisfies
$$
\lim_{t\to t_0^+} (u(t,\cdot;t_0,u_0,v_0),v(t,\cdot;t_0,u_0,v_0))=(u_0(\cdot),v_0(\cdot))
$$
{ in $C(\bar\Omega)\times W^{1,q}(\Omega)$ for any $q>n.$ } In such case, $(u(t,x;t_0,u_0,v_0),v(t,x;t_0,u_0,v_0))$ is called the solution
of \eqref{u-v-eq00} on $[t_0,t_0+T)$ with initial condition $(u(t_0,x),v(t_0,x))=(u_0(x),v_0(x))$. If $T=\infty$, we say
the solution of \eqref{u-v-eq00} with initial condition $(u_0(x),v_0(x))$ at $t=t_0$ exists globally.
\end{definition}

Next, we definite the notions of the entire solution and pointwise persistence.
\begin{definition}
\label{coex-persist-def}
A solution $(u(t,x),v(t,x))$ of \eqref{u-v-eq00} defined for all $t\in\RR$ is called an {\rm entire solution}.
A positive entire solution of \eqref{u-v-eq00} is an entire  solution $(u^{**}(t,x),v^{**}(t,x))$ satisfying \begin{equation*}
\label{coexistence-eq}
\inf_{t\in\RR,x\in\bar\Omega} u^{**}(t,x)>0.
\end{equation*}
We say that {\rm pointwise persistence} occurs in \eqref{u-v-eq00} if there is $\eta>0$ such that for any $t_0\in\RR$, $u_0\in C(\bar\Omega),$ and $v_0\in W^{1,\infty}(\Omega)$ with $u_0> 0$ and  $v_0\geq 0$, $(u(t,x;t_0,u_0,v_0),v(t,x;t_0,u_0,v_0))$ exists globally, and
there is $\xi(t_0,u_0,v_0)>0$ such that
\begin{equation*}
\label{persistence-eq}
u(t,x;t_0,u_0,v_0)\ge \eta,\quad \forall\,\, t\ge t_0+\xi(t_0,u_0,v_0).
\end{equation*}
\end{definition}

{Thanks to the results established in \cite{ITBWS2020}, throughout the rest of this work, we shall always assume that the following  hypothesis  holds:
\medskip

{\bf (H1)} \quad$ \inf_{t\in\mathbb{R}}(a_{1,\inf}(t)-|\Omega|{(a_{2,\inf}(t))_{-}})>0.$ 
}

\medskip

\noindent Under hypothesis {\bf (H1)}, the following result on the global boundedness, persistence of positive classical solutions, and existence of positive entire solutions of system \eqref{u-v-eq00} is established in \cite{ITBWS2020}.

\begin{theorem}\cite[Theorems 1.1, 1.2 and 1.3]{ITBWS2020}
\label{thm-global-000}
  Assume that {\bf (H1)} holds.  Then there exists $\chi_0={\chi_0(a_i,\tau,\Omega,\mu,\lambda)}>0$, $M=M(\chi_0)>0$ and $\eta=\eta(\chi_0)>0$ such for every $|\chi|\leq \chi_0$, $t_0 \in \R,$ $(u_0,v_0) \in C(\bar\Omega )\times W^{1,\infty}(\Omega)$ { with $u_0,v_0 > 0$}, \eqref{u-v-eq00} has a unique bounded {globally defined} classical solution $(u(t,\cdot;t_0,u_0,v_0)$, $v(t,\cdot;t_0,u_0,v_0))$.  Furthermore, there exists $t^1(u_0,v_0)>0$ such that
  \begin{equation}
  \label{m2-eq}
  \inf_{x \in \bar{\Omega}}u(t,x;t_0,u_0,v_0) \geq \eta \quad \text{ and } \quad \|u(t,\cdot;t_0,u_0,v_0)\|_{\infty}\le M\quad \forall\, t\ge t_0+t^1(u_0,v_0).
  \end{equation}
Moreover, system \eqref{u-v-eq00} has at least one positive entire solution ${(u^*(t,x),v^*(t,x))}$.

\end{theorem}
{
\begin{proof} We take 
$ \chi_0:=\frac{a_{1,\inf}}{\inf_{q>\max\{1,\frac{n}{2}\}} \Big(\frac{q-1}{q}(C_{q+1})^{\frac{1}{q+1}}\mu^{\frac{1}{q+1}}\Big)}
$, where the positive constants $C_{q+1}$ are given in \cite[Lemma 1.1]{ITBWS2020}. Hence $\chi_0$ together with our hypothesis {\bf (H1)} in the crrent work satisfies \cite[Hypothesis {\bf (H1)}]{ITBWS2020}. Now the results, follows from \cite[Theorems 1.1, 1.2 and 1.3]{ITBWS2020}
    
\end{proof}
}
\subsection{Main result}
We state the main result of the manuscript in the current subsection. Let $\eta$ and $M$ as in Theorem \ref{thm-global-000}. To establish the uniqueness and global stability of positive entire solutions of system \eqref{u-v-eq00}, the following additional hypothesis is needed:

\medskip
\noindent {\bf (H2)}  { $\eta a_{1,\inf} >|\Omega|M((a_2)_{+,\sup}+(a_2)_{-,\sup}).$}

\medskip

Our main result on the uniqueness and global stability of the positive entire solutions of system \eqref{u-v-eq00} reads as follows:

\begin{theorem}\label{thm-nonlinear-stability}{ (Uniqueness and stability of entire solutions)}
Assume  that {\bf (H1)} and {\bf (H2)} { hold}. Then there is $\chi_1>0$ such that for every $|\chi|\leq \chi_1,$ system \eqref{u-v-eq00} has a unique {positive entire solution} $(u^*(t,x),v^*(t,x))$ which is globally stable in the sense that  
  \begin{equation}
  \label{global-stability-2-eq1}
  \lim_{t \to \infty}\Big(\sup_{t_0 \in \mathbb{R}}\|u(t{+t_0},\cdot;t_0,u_0,v_0)-u^*(t{+t_0},\cdot)\|_{\infty}+\sup_{t_0 \in \mathbb{R}}\|v(t{+t_0},\cdot;t_0,u_0,v_0)-v^*(t{+t_0},\cdot)\|_{\infty}\Big)=0,
  \end{equation}
for any  $u_0 \in C^0(\bar{\Omega}),$ $v_0 \in W^{1,\infty}(\bar{\Omega})$ with $u_0,v_0\ge 0$ and
$u_0\not\equiv 0$.

\end{theorem}




}

We conclude this section with the following remarks. 

\begin{remark}

Theorem \ref{thm-nonlinear-stability} extends the uniqueness and stability of entire
solutions result obtained by Winkler in  \cite{Win2014} for a parabolic-parabolic chemotaxis model with { $\tau=1$ } in homogeneous convex environments to nonconvex bounded heterogeneous environments with {$\tau>0$}. This extension is fairly {nontrivial} and is done mostly by applying a {nontrivial} application of the so-called method of eventual comparison adapted to heterogeneous environments.

\end{remark}

\medskip

The rest of the paper is organized as follows. In section \ref{prem},  we recall crucial lemmas, and in section \ref{uniquess}, we will present the proof of Theorem \ref{thm-nonlinear-stability}.

\section{Preliminary}
\label{prem}
Throughout the rest of this paper, we shall suppose that {\bf (H1)} and {\bf (H2)} hold and the positive constants $\chi_0,\eta$ and $M$ are as in Theorem \ref{thm-global-000}. We start with the following lemma on boundedness of $v$ in $W^{2,\infty}(\Omega)$ eventually as time becomes large.

\begin{lemma}\label{lem-00}
There exist $C_0=C(\chi_0)$ {and $\zeta^*=\zeta^*(\chi_0)>0$}  such that if $|\chi|\leq \chi_0,$
then for any  $t_0 \in \mathbb{R}$, and $(u_0,v_0) \in C(\bar{\Omega})\times W^{1,\infty}(\Omega)$ with $u_0,v_0 \geq,\ne 0,$  there is $t^*=t^*(u_0,v_0)>0$ satisfying that
\begin{equation}
\label{new-new-eq3-1}
\|v(\cdot,t)\|_{W^{2,\infty}(\Omega)}\le C_0 \sup_{t_0+t^*\le s\le t}\|u(\cdot,s)\|_\infty\quad \forall \, t\ge t_0+t^*{ +\zeta^*}.
\end{equation}
\end{lemma}
\begin{proof}

   \noindent{ Fix $p>n$, $0<\alpha<1$, $0<\zeta<\frac{1}{2}$, and set $Aw:=\alpha w-\Delta w$ for all $w\in D(\Delta_p):=\{w\in W^{2,p}(\Omega) : \partial_{\nu}w=0\ \text{on }\ \partial\Omega\}$. Let $A^{\zeta}$ denote the fractional power of $A$.  By  \eqref{m2-eq},  there is $t^1(u_0,v_0)\gg 0$ such that 
    \begin{equation}\label{RA1}
       \eta\le u(t,x;t_0,u_0,v_0)\le M \quad \forall\ t\ge t_0+t^1(u_0,v_0).
    \end{equation}
    Next, by \cite[Lemma 4.3]{ITBWS2020}, there exist  $C^*_1>0$ and  $\theta^*>0$, independent of initial data such that 
    \begin{equation}\label{RA2}
        \|v(\cdot,t) \|_{W^{2,\infty}(\Omega)}\le C^*_1(t-\tilde{t}_0)^{-\theta^*}e^{\frac{-(\lambda{-\alpha)}}{\tau}(t-\tilde{t}_0)}\|v(\cdot,\tilde{t}_0)\|_{\infty}+C^*_1\sup_{\tilde{t}_0\le s\le t}\|A^{\zeta}u(\cdot,s)\|_{L^p(\Omega)}\quad \forall\ t>\tilde{t}_0>t_0.
    \end{equation}
    By \cite[Lemma 4.2]{ITBWS2020}, there is $C^*_2>0$, independent of initial data, such that 
    \begin{align}\label{RA3}
&\|A^{\zeta}u(\cdot,t)\|_{L^p(\Omega)}\le  C^*_2(t-\tilde{t}_0)^{-\zeta}e^{-(1-\alpha)(t-\tilde{t}_0)}\|u(\cdot,\tilde{t}_0)\|_{L^p(\Omega)}\cr
        &+ C^*_2\Big(1+\sup_{\tilde{t}_0\le s\le t}\|\nabla v(\cdot,s)\|_{L^p(\Omega)}+\sup_{\tilde{t}_0\le s\le t}\|u(\cdot,s)\|_{\infty}\Big)\sup_{\tilde{t}_0\le s\le t}\|u(\cdot,s)\|_{\infty}\quad\forall\ t>\tilde{t}_0>t_0.
    \end{align}
    By \cite[Corollary 4.3]{ITBWS2020}, there is  $C_3^*>0$, independent of initial data, and $t^2(u_0,v_0)>0$ such that 
    \begin{equation}\label{RA4}
        \|v(\cdot,t)\|_{W^{2,\infty}(\Omega)}\le C_3^*\quad \forall\ t\ge t_0+t^2(u_0,v_0).
    \end{equation} %
    Taking $C_0^*:=\max\{C_1^*,C_2^*,C_3^*\}(1+(n|\Omega|)^{\frac{1}{p}}C_3^*+M)$, $\tilde{\zeta}_*:=\min\{\zeta,\theta^*\}$, $\tilde{\lambda}_*=\min\{(1-\alpha),\frac{(\lambda{-\alpha)}}{\tau}\}$, and   $t^3(u_0,v_0)=t^1(u_0,v_0)+t^2(u_0,v_0)$, it follows from \eqref{RA1}, \eqref{RA3}, and \eqref{RA4} that 
    \begin{equation*}
        \|A^{\zeta}u(\cdot,t)\|_{L^p(\Omega)}\le C_0^*\left(|\Omega|^{\frac{1}{p}}M(t-\tilde{t}_0)^{-{\tilde{\zeta}^*}}e^{-\tilde{\lambda}_*(t-\tilde{t}_0)}+\sup_{\tilde{t}_0\le s\le t}\|u(\cdot,s)\|_{\infty}\right)\quad \forall\ t\ge 1+\tilde{t}_0,
    \end{equation*}
    for any $\tilde{t_0}\ge t_0+t^3(u_0,v_0)$. This along with \eqref{RA2} and \eqref{RA4} gives that,  for any $\tilde{t}_0\ge t_0+t^3(u_0,v_0)$,
    \begin{equation}\label{RA5}
        \|v(\cdot,t)\|_{W^{2,\infty}(\Omega)}\le C_0^*\left({(}C_3^*+|\Omega|^{\frac{1}{p}}MC_0^*)e^{-\tilde\lambda_*(t-\tilde{t}_0)}+C_0^*\sup_{\tilde{t}_0\le s\le t}\|u(\cdot,s)\|_{\infty}\right)\quad \forall\ t\ge 1+\tilde{t}_0.
    \end{equation}
   Therefore, choosing $\zeta_0\gg 1$ such that
    $$
    (C_3^*+|\Omega|^{\frac{1}{p}}MC_0^*)e^{-\tilde\lambda_*\zeta_0}<C_0^*\eta,
    $$ we conclude from \eqref{RA1} and \eqref{RA5} that 
    \begin{align*}
    \|v(\cdot,t)\|_{W^{2,\infty}(\Omega)}
\le2[C^*_0]^2\sup_{{\tilde{t}_0}\le s\le t}\|u(\cdot,s)\|_{\infty}\quad { \forall\ t\ge 1+\zeta_0+\tilde{t}_0 } .
    \end{align*}
    Therefore \eqref{new-new-eq3-1} holds with $C_0=2[C_0^*]^2$, {$\zeta^*=1+\zeta_0$, and $t^*(u_0,v_0)=
    t^3(u_0,v_0)$.}
    }
\end{proof}

Next, we state a result that gives a uniform bound of gradient of $u^*$ whenever $(u^*,v^*)$ is a positive entire solution of system \eqref{u-v-eq00}. First,  for every $p>1$ and $\frac{1}{2}<\beta<1$, let $X^{\beta}$ denote the fractional power space induced by the invertible, closed and densely defined linear operator ${\rm id}-\Delta$ on $L^p(\Omega)$. In the following, we fix $p>n$ and $\frac{1}{2}<\beta<1$ such that $1+\frac{n}{p}<2\beta$. Then by \cite[Theorem 1.6.1]{Dan}, $X^{\beta}$ is compactly embedded in $C^1(\bar{\Omega})$. 

\begin{lemma}
\label{lem-01}
Let  $\frac{1}{2}<\beta<1$ and $p>n$ such that $1+\frac{n}{p}<2\beta$. Then there exists a positive constant $C_1=C_1(\chi_0)$ { such that} if {$|\chi|\le \chi_0$ and}  $(u^*(t,x),v^*(t,x){ )}$ is { a positive} entire solution of system \eqref{u-v-eq00}, we have:
\begin{equation}
\label{new-new-eq3-1}
    \|u^*(t)\|_{X^{\beta}}\le C_1,\quad \forall\ t\in  \mathbb{R}.
\end{equation}
In particular, since $X^{\beta} \xhookrightarrow{} C^1(\bar{\Omega}),$ by equation \eqref{new-new-eq3-1}, there exists $\tilde C_1=\tilde C_1(\chi_0)$ such that
\begin{equation}
\label{new-new-eq3-2}
    \|u^*(t)\|_{C^1(\bar{\Omega})}\le \tilde C_1, \quad \forall\ t\in  \mathbb{R}. 
\end{equation}
\end{lemma}
\begin{proof}
  {
Let $(u^*(t,x),v^*(t,x){ )}$ be a positive entire solution of system \eqref{u-v-eq00}. Then $u^*(t,x)$ { satisfies} the following equation
$$
\begin{cases}
u_t^*=(\Delta -1) u^* +\nabla u^*\cdot\underbrace{\nabla (-\chi v^*)}_{B}+\underbrace{(1+a_0-a_1u^*-a_2\int_{\Omega}u^*-\chi\Delta v^*)u^*}_{A}  & x\in\Omega,\ t\in\mathbb{R},\cr
\frac{\partial u^*}{\partial\nu}=0 & x\in\partial\Omega, \ t\in\mathbb{R},
\end{cases}
$$
and by the variation of constants formula
\begin{equation}\label{FGR1-1}
u^*(t+t_0)=e^{t(\Delta -1)}u^*(t_0)+\int_{0}^te^{(t-s)(\Delta -1)}(\nabla u^*(t_0+s)\cdot B(t_0+s)+A(t_0+s))ds\quad \forall\ t\in\mathbb{R}.
\end{equation}
Next, thanks to \cite[Theorem 1.4.3]{Dan}, for every $0<\delta<1$, there exists a positive constant $K_{1}=K_1(\beta,\delta,p,n)>0 $ such that 
$$
\|e^{t(\Delta-1)}z\|_{X^{\beta}}\le K_1t^{-\beta}e^{-\delta t}\|z\|_{L^p(\Omega)}\quad \forall\ t>0.
$$
Since $X^{\beta}$ is continuously embedded in $C^1(\bar{\Omega})$ (\cite[Theorem 1.6.1]{Dan}), there is $K_2>0$ such that 
$$
\|z\|_{C^{1}(\bar{\Omega})}\le K_2\|z\|_{X^{\beta}} \quad \forall\ z\in X^{\beta}.
$$
Fix $0<\delta<1$. Then, by \eqref{FGR1-1},  for every $0<t\le 1$ and $t_0\in\mathbb{R}$, we have
\begin{align*}
&\|u^*(t+t_0)\|_{X^{\beta}}\cr 
\le& \|e^{t(\Delta -1)}u^*(t_0)\|_{X^{\beta}}+\int_0^t\|e^{(t-s)(\Delta -1)}(\nabla u^*(t_0+s)\cdot B(t_0+s)+A(t_0+s))\|_{X^{\beta}}ds\cr
\le & \frac{|\Omega|^{\frac{1}{p}}K_1}{t^{\beta}e^{t\delta}}\|u^*\|_{\sup}+K_1\int_0^t\frac{e^{-\delta(t-s)}}{(t-s)^\beta}(\|B(t_0+s)\|_{\infty}\||\nabla u^*(t_0+s)|\|_{L^p(\Omega)}+\|A(t_0+s)\|_{L^p(\Omega)})ds\cr 
\le & \frac{K_1|\Omega|^{\frac{1}{p}}}{t^{\beta}e^{\delta t}}\|u^*\|_{\sup}+K_1|\Omega|^{\frac{1}{p}}\|A\|_{\sup}\int_{0}^t\frac{e^{-\delta(t-s)}}{(t-s)^{\beta}}ds+n^{\frac{1}{p}}|\Omega|^{\frac{1}{p}}K_1\|B\|_{\sup}\int_0^t\frac{e^{-\delta(t-s)}}{(t-s)^{\beta}}\|u^*(t_0+s)\|_{C^1(\bar{\Omega})}ds\cr 
\le & \frac{K_1|\Omega|^{\frac{1}{p}}}{t^{\beta}e^{\delta t}}\|u^*\|_{\sup}+K_1|\Omega|^{\frac{1}{p}}\|A\|_{\sup}\delta^{\beta-1}\Gamma(1-\beta)+n^{\frac{1}{p}}|\Omega|^{\frac{1}{p}}K_1\|B\|_{\sup}\int_0^t\frac{e^{-\delta(t-s)}}{(t-s)^{\beta}}\|u^*(t_0+s)\|_{C^1(\bar{\Omega})}ds\cr 
\le& K_1|\Omega|^{\frac{1}{p}}\Big({\|u^*\|_{\sup}}+\|A\|_{\sup}\delta^{\beta-1}\Gamma(1-\beta)\Big)\frac{1}{t^{\beta}}+K_2n^{\frac{1}{p}}|\Omega|^{\frac{1}{p}}K_1\|B\|_{\sup}\int_0^t\frac{1}{(t-s)^{\beta}}\|u^*(t_0+s)\|_{X^{\beta}}ds.
\end{align*}
where  
$\Gamma(1-\beta)$ is the Gamma function evaluated at $1-\beta$. Note from Theorem \ref{thm-global-000} and Lemma \ref{lem-00} that there is $M_1=M_{1}(\chi_0)>0$ such that 
$$
\|u^*\|_{\sup}\le M_1, \quad \|B\|_{\sup}\le |\chi|M_1 \quad \text{and}\quad \|A\|_{\sup}\le M_1 \quad \forall\ |\chi|\le \chi_0. 
$$
This along with the previous inequality gives
$$
\|u^*(t_0+t)\|_{X^{\beta}}\le M_2t^{-\beta}+M_3\int_{0}^t\frac{\|u^*(t_0+s)\|_{X^{\beta}}}{(t-s)^\beta}ds\quad \forall\ 0<t\le 1,\ |\chi|\le \chi_0,\ t_0\in\mathbb{R}.
$$
where $M_2:=K_1M_1|\Omega|^{\frac{1}{p}}(1+\delta^{\beta-1}\Gamma(1-\beta))>0$ and $M_3:=\chi_0n^{\frac{1}{p}}M_1K_1K_2|\Omega|^{\frac{1}{p}}>0$. Thus, by \cite[page 190-Exercise $4^*$]{Dan}, there { exists} a positive number $K_3=K_3(M_2,M_3,\beta)>0$ such that 
$$
\|u^*(t+t_0)\|_{X^{\beta}}\le K_3t^{-\beta}\quad \forall\ 0<t<1,\ \   |\chi|\le \chi_0, \ t_0\in\mathbb{R}.
$$
Therefore, \eqref{new-new-eq3-1} holds with $C_1=K_3$.
}
 
\end{proof}

\begin{lemma}
\label{lem-02}
Let $(u^*(t,x),v^*(t,x))$ { be} a positive entire solution of system \eqref{u-v-eq00}. 
There { exist} $C_2=C_2(\chi_0)$  and $\tilde{\lambda}_0>0$
such that if $|\chi|\leq \chi_0,$ then
for any  $t_0 \in \mathbb{R}$, and $(u_0,v_0) \in C(\bar{\Omega})\times W^{1,\infty}(\Omega)$ with $u_0,v_0 \geq 0,$  there is $t^*=t^*(u_0,v_0)>0$ satisfying that
\begin{align}
\label{delta-bound-eq0}
  &|\Delta(v^*-v)(t)+\nabla \ln{u^*}(t)\cdot\nabla(v^*-v)(t)|\cr
  \leq & C_2\Big(e^{-\tilde{\lambda}_0(t-t_0-t^*(u_0,v_0))} +\sup_{t_0+t^*\le s\le t}\|u(s,\cdot){-u^*(s,\cdot)}\|_\infty\Big)\quad \forall \, t\ge 1+ t_0+t^*.  
\end{align}

\end{lemma}
\begin{proof}
Let $t \in \mathbb{R}$, then 
\begin{align*}
  &|\Delta(v^*-v)(t)+\nabla \ln{u^*}(t)\cdot\nabla(v^*-v)(t)| \leq \Big(1+|\nabla \ln{u^*}(t)|\Big){ \|(v^*-v)(t)\|_{W^{2,\infty}(\Omega)} }\cr 
  =&\Big(1+\frac{|\nabla u^*(t)|}{|u^*(t)|}\Big){ \|(v^*-v)(t)\|_{W^{2,\infty}(\Omega)}}\leq \Big(1+\frac{\tilde C_1}{\eta}\Big){ \|(v^*-v)(t)\|_{W^{2,\infty}(\Omega)}},
\end{align*}
where $\tilde C_1$ is given by Lemma \ref{lem-01}. The  desired result follows  as in \eqref{RA5}. 
\end{proof}

\section{Uniqueness and stability of positive entire solutions}
\label{uniquess}
In this section,  we established the proof of Theorem \ref{thm-nonlinear-stability}.{ We  prove three crucial lemmas. 

\begin{lemma}
    \label{LL0}
    Let $(u^*(t,x),v^*(t,x))$ be a positive entire solution of system \eqref{u-v-eq00} and given $t_0\in\RR$, $u_0\in C(\bar\Omega)$, and $v_0\in W^{1,\infty}(\Omega)$ with $u_0> 0$ and  $v_0\geq 0$, let $(u(t,x),v(t,x))=(u(t,x;t_0,u_0,v_0),v(t,x;t_0,u_0,v_0))$. Then if $w=\frac{u}{u^*}$ and $\phi=\frac{v}{v^*},$ $(w,\phi)$ satisfies the following system of equations:
\begin{equation}
 \label{w-phi-eq00}
\begin{cases}
w_t=\Delta w+\nabla w \cdot\nabla [2 \ln{u^*}-\chi v ]+\chi[\Delta(v^*-v)+\nabla \ln{u^*}\cdot\nabla(v^*-v)]w- F(a_1,a_2,w,u^*),& x\in \Omega,\cr
\tau \phi_t=\Delta \phi+2\nabla \phi\cdot \nabla \ln{v^*}+\mu\frac{u^*}{v^*}(w-\phi),& x\in \Omega, \cr
\frac{\p w}{\p \nu}=\frac{\p \phi}{\p \nu}=0,& x\in\p\Omega,
\end{cases}
 \end{equation}
 where ${ F(a_1,a_2,w,u^*)}=\big(a_1(t,x)u^*(w-1)+a_2(t,x)\int_{\Omega}u^*(w-1)\big)w.$

\end{lemma}
\begin{proof}
Note that $u^*w=u,$ so $u^*w_t=u_t-u^*_tw$. Since both $u$ and $u^*$ satisfy the first equation of system \eqref{u-v-eq00}, we get
\begin{align}
\label{w-eq0}
u^*w_t&=\Delta u-\chi\nabla\cdot (u \nabla v)+\big(a_0(t,x)-a_1(t,x)u-a_2(t,x)\int_{\Omega}u\big)u\\ \nonumber
&-\big[\Delta u^*-\chi\nabla\cdot (u^* \nabla v^*)+\big(a_0(t,x)-a_1(t,x)u^*-a_2(t,x)\int_{\Omega}u^*\big)u^*\big]w
\end{align}
Next dividing by $u^*$ on { both sides} of equation \eqref{w-eq0}, we get
\begin{align}
\label{w-eq1}
w_t=\frac{\Delta u-w\Delta u^*}{u^*}-\chi\frac{\nabla\cdot (u \nabla v)}{u^*}+\chi\frac{w\nabla\cdot (u^* \nabla v^*)}{u^*}-\big(a_1(t,x)(u-u^*)+a_2(t,x)\int_{\Omega}(u-u^*)\big)w
\end{align}
Note that $\frac{\Delta u-w\Delta u^*}{u^*}=\Delta w+2 \nabla w \cdot \nabla \ln{u^*}$, 
\begin{align*}
  \frac{\nabla\cdot (u \nabla v)}{u^*} =  \frac{\nabla u \cdot \nabla v}{u^*} +w \Delta v 
  &=w \frac{\nabla u^* \cdot \nabla v}{u^*}+\nabla w \cdot \nabla v +w \Delta v\\ \nonumber
  &=w \nabla \ln{u^*} \cdot \nabla v+\nabla w \cdot \nabla v +w \Delta v,
\end{align*}
and
\begin{align*}
  \frac{\nabla\cdot (u^* \nabla v^*)}{u^*}& =  \frac{\nabla u^* \cdot \nabla v^*}{u^*} + \Delta v^*=\nabla \ln{u^*} \cdot \nabla v^*+\Delta v^*.
\end{align*}
Equation \eqref{w-eq1} becomes
\begin{align*}
w_t&=\Delta w+ \nabla w \cdot \nabla [2\ln{u^*}-\chi v]+\chi[\Delta(v^*-v)+\nabla \ln{u^*} \cdot \nabla(v^*-v)]w\\ \nonumber
&+\big(a_1(t,x)u^*(1-w)+a_2(t,x)\int_{\Omega}u^*(1-w)\big)w.
\end{align*}
Similarly, since both $v$ and $v^*$ satisfy the second equation of system \eqref{u-v-eq00}, { we get}
\begin{align*}
    \tau \phi_t &=\frac{\Delta v-\lambda v +\mu u}{v^*}-\frac{\Delta v^*-\lambda v^* +\mu u^*}{v^*}\phi \\ \nonumber
    &=\frac{\Delta v-\phi  \Delta v^* }{v^*}+\mu\frac{u-u^*\phi}{v^*} -\lambda (\frac{v}{v^*}-\phi)= \Delta \phi +2 \nabla \phi \cdot \nabla \ln{v^*}+  { \mu}\frac{u^*}{v^*}(w-\phi).
  \end{align*}
\end{proof}

Fix $t_0 \in \mathbb{R}$, and $(u_0,v_0) \in C(\bar{\Omega})\times W^{1,\infty}(\Omega)$ with $u_0,v_0 > 0.$ Then by Theorem \ref{thm-global-000} and Lemma \ref{lem-02}, there exists $T^*>0$ independent of $t_0$ { such that}
\begin{equation}
\label{gradients-bounds}
  |\Delta(v^*-v)(t{+t_0},\cdot)+\nabla \ln{u^*}(t{+t_0},\cdot)\cdot\nabla(v^*-v)(t{+t_0},\cdot)| \leq 2C_2M:=C\quad \forall \, t\ge T^*,\ {t_0\in\mathbb{R}},   
\end{equation}
where $C_2$ is given by Lemma \ref{lem-02} and $M$ by Theorem \ref{thm-global-000}.

To establish an eventual pointwise estimate on $w(t,x)$, we consider the system of two species competition system  

 \begin{equation}\label{ER0}
     \begin{cases}
         \frac{d\overline{w}}{dt}=\Big[|\chi|C+ \underbrace{\eta a_{1,\inf}}_{A_0}(1-\overline{w})+\underbrace{|\Omega|M(a_2)_{-,\sup}}_{A_{1}}(\overline{w}-1)+\underbrace{|\Omega|M(a_2)_{+,\sup}}_{A_{2}}(1-\underline{w})\Big]\bar{w} & t>T^*\cr
          \frac{d\underline{w}}{dt}=\Big[-|\chi|C+ \underbrace{\eta a_{1,\inf}}_{A_0}(1-\underline{w})+\underbrace{M|\Omega|(a_2)_{-,\sup}}_{A_{1}}(\underline{w}-1)+\underbrace{M|\Omega| (a_2)_{+,\sup}}_{A_{2}}(1-\bar{w})\Big]\underline{w}& t>T^*\cr
          (\bar{w}(T^*),\underline{w}(T^*))=\big(w_{\sup}(T^*,\cdot)+1,\frac{w_{\inf}(T^*,\cdot)}{2w_{\inf}(T^*,\cdot)+1}\big),
     \end{cases}
 \end{equation}
where $C$ is given by equation \eqref{gradients-bounds} and $\eta, M$ by Theorem \ref{thm-global-000}, and { $w_{\sup}(T^*,\cdot):=\sup_{t_0\in\mathbb{R}}\|w(T^*+t_0,\cdot)\|_{\infty}$ and $w_{\inf}(T^*,\cdot):=\inf_{t_0\in\mathbb{R}}\min_{x\in\bar{\Omega}}w(T^*+t_0,x)$}.  Note that $w_{\inf}(T^*,\cdot)>0$ by Theorem \ref{thm-global-000} and system \ref{ER0} is equivalent to 
  \begin{equation}\label{ER0-bis}
     \begin{cases}
         \frac{d\overline{w}}{dt}=\Big[|\chi|C+ A_0-A_{1}+A_{2}-(A_0-A_1)\overline{w}-A_{2}\underline{w}\Big]\overline{w} & t>T^*\cr
          \frac{d\underline{w}}{dt}=\Big[-|\chi|C+  A_0-A_{1}+A_{2}-(A_0-A_{1})\underline{w}-A_{2}\overline{w}\Big]\underline{w}& t>T^*\cr
          (\bar{w}(T^*),\underline{w}(T^*))=\big(w_{\sup}(T^*,\cdot)+1,\frac{w_{\inf}(T^*,\cdot)}{2w_{\inf}(T^*,\cdot)+1}\big),
     \end{cases}
 \end{equation}

 It is clear from the comparison principle for two species competition systems that $\underline{w}(t)<1<\bar{w}(t)$ for all $t\ge T^*$.  The following result on the large time behavior of solutions of \eqref{ER0} follows from classical results on the dynamics of two species Lotka-Voltera competition systems (see for example \cite{SAhmad1987,ITBWS17b,JMS1968}). Note that { {\bf (H2)} holds if and only if  $A_0-A_1-A_2>0.$}

\begin{lemma}\label{LL1-1} Let $A_{i}$, $i=0,1,2$ be as in \eqref{ER0}. Suppose that {\bf (H1)} and {\bf (H2)} hold. Let $|\chi|<\min\{\frac{A_0-A_1-A_2}{{ 2}C},\chi_0\}$. Then, the solution $(\bar{w}(t), \underline{w}(t))$ of \eqref{ER0} converges to $(\bar{w}^*,\underline{w}^*)$ as $t\to\infty$, where 
\begin{small}
\begin{equation}
    \overline{w}^*=1+\frac{C|\chi|}{A_0-A_{1}-A_{2}}>0\quad \text{and}\quad \underline{w}^*=1-\frac{C|\chi|}{A_0-A_{1}-A_{2}}>0.
\end{equation}
\end{small}
    
\end{lemma}

Our next result gives a lower and upper bounds for $w(t,x)$.

\begin{lemma}\label{LL2} Let $(\bar{w}(t),\underline{w}(t))$ be the unique solution of \eqref{ER0}. Then 
\begin{equation}\label{LL2-eq1}
    \underline{w}(t)<w(t{+t_0},x)<\bar{w}(t)\quad \forall\ t\ge T^*,\ {t_0\in\mathbb{R}},\ x\in\bar{\Omega}. 
\end{equation}
    
\end{lemma}
 \begin{proof} Suppose to the contrary that \eqref{LL2-eq1} is false. Then there exit {$t_0\in\mathbb{R}$ and} $T^{**}>T^*$ such that \eqref{LL2-eq1} holds $[T^*, T^{**})$ and either 
 \begin{equation}\label{ER7}
     w(T^{**}{+t_0},x_{**})=\overline{w}(T^{**})\quad  \text{or} \quad w(T^{**}{+t_0},x_{**})=\underline{w}(T^{**}) \ \text{for some}\ x_{**}\in\bar{\Omega}.
 \end{equation}
Now, since $\underline{w}<1<\overline{w}$ and $\underline{w}(t)\le { w(t+{t_0},\cdot)}\le \overline{w}(t)$ for $t\in [T^*,T^{**}]$, setting $B:=\nabla [2 \ln{u^*}-\chi v]$,  we have  
\begin{align}
 \label{w-upper-equation-0}
w_t(t{+t_0},\cdot) \leq& \Delta w+B\cdot\nabla w + |\chi|Cw + \big(a_1(t,x)u^*(1-w)+a_2(t,x)\int_{\Omega}u^*(1-w)\big)w\cr 
=& \Delta w +B\cdot\nabla w +\Big(|\chi|C +a_1u^*(1-w)+(a_2)_+\int_{\Omega}u^*(1-w)+(a_2)_{-}\int_{\Omega}u^*(w-1)\Big)w\cr
\le &\Delta w +B\cdot\nabla w +\Big(|\chi|C +a_1u^*(1-{w})+\Big((a_2)_+\int_{\Omega}u^*\Big)(1-\underline{w})+\Big((a_2)_{-}\int_{\Omega}u^*\Big)(\overline{w}-1)\Big)w\cr
\le &\Delta w +B\cdot\nabla w +\Big(|\chi|C +a_1u^*(1-{w})+|\Omega|M(a_2)_{+,\sup}(1-\underline{w})+|\Omega|M(a_2)_{-,\sup}(\overline{w}-1)\Big)w\cr
=&\Delta w +B\cdot\nabla w +\Big(|\chi|C +a_1u^*(1-{w})+A_{1}(\overline{w}-1)+A_{2}(1-\underline{w})\Big)w,\quad T^*\le t\le T^{**}.\cr 
 \end{align}
 Similarly, 
 it holds that 
\begin{align}
 \label{w-upper-equation-1}
w_t(t{+t_0},\cdot) \geq&\Delta w+B\cdot\nabla w +\Big(-|\chi|C + a_1(t,x)u^*(1-w)+(a_2)_+\int_{\Omega}u^*(1-\bar{w})+(a_{2})_{-}\int_{\Omega}u^*(\underline{w}-1)\Big)w\cr 
\ge & \Delta w +B\cdot { \nabla} w +\Big(-|\chi|C+a_1u^*(1-w) +M|\Omega|(a_2)_{+,\sup}(1-\overline{w})+M|\Omega|(a_2)_{-,\sup}(\underline{w}-1)\Big)w\cr
=&\Delta w +B\cdot { \nabla} w +\Big(-|\chi|C+a_1u^*(1-w) +A_{2}(1-\overline{w})+A_{1}(\underline{w}-1)\Big)w, \quad T^*\le t\le T^{**}.\quad 
 \end{align}
 Next, recalling that $\overline{w}(t)>1$ for all $t\ge T^*$, it follows from the first equation of \eqref{ER0} that 
 $$
 \frac{d\overline{w}}{dt}\ge \Big(|\chi|C+a_1u^*(1-\overline{w})+A_{1}(\overline{w}-1)+A_{2}(1-\underline{w})\Big)\bar{w}\quad t>T^*. 
 $$
 This along with \eqref{w-upper-equation-0}, the facts that $\partial_{\nu}w(t{+t_0})=0=\partial_{\nu}\overline{w}$ for $t\in (T^*,T^*]\times \bar{\Omega}$, $w(T^*+{ t_0},\cdot)<\overline{w}(T^*)$, and the strong maximum principle for { single parabolic equation} give 
 \begin{equation*}
     w(t{+t_0},x)<\overline{w}(t)\quad \forall\ T^*<t\le T^{**}, \ x\in\bar{\Omega}.
 \end{equation*}
 In particular $w(T^{**}{+t_0},x_{**})<\bar{w}(T^{**})$, shows that the first possibility in \eqref{ER7} cannot hold. Similarly, recalling that ${ \underline{w}(t)<1}$ for all $t>T^*$, 
 it follows from the second equation of \eqref{ER0} that 
 $$
 \frac{d\underline{w}}{dt}\le \Big(-|\chi|C+a_1u^*(1-\underline{w})+A_{1}(\underline{w}-1)+A_{2}(1-\overline{w})\Big)\underline{w}\quad t>T^*. 
 $$
 This along with \eqref{w-upper-equation-1}, the facts that $\partial_{\nu}w(t{+t_0},\cdot)=0=\partial_{\nu}\underline{w}$ on $(T^*,T^*]\times \bar{\Omega}$, $w(T^*+{ t_0},\cdot)>\underline{w}(T^*)$, and the strong maximum principle for { single parabolic equation} give 
 \begin{equation*}
     w(t{+t_0},x)>\underline{w}(t)\quad \forall\ T^*<t\le T^{**}, \ x\in\bar{\Omega}.
 \end{equation*}
 In particular $w(T^{**}{+t_0},x_{**})>\underline{w}(T^{**})$, which shows that the second possibility in \eqref{ER7} cannot also hold. Therefore, we have that \eqref{ER7} cannot hold, and hence we must have that $T^{**}=\infty$, which completes the proof of the lemma.

\end{proof}

Finally, we give a proof of  Theorem \ref{thm-nonlinear-stability}.
\begin{proof}
Let $|\chi|<\chi_1:=\min\{\frac{A_0-A_1-A_2}{C},\chi_0\}$. Then given $0<\epsilon<1$ and $n \in \mathbb{N}$, we claim that there exists $T_{n,\epsilon}{\ge n}$ such that
\begin{equation}
\label{uniform-conv-eq00}
    \sup_{t_0 \in \mathbb{R}}\|u(t{+t_0},\cdot)-u^*(t{+t_0},\cdot)\|_\infty\leq M\frac{(C|\chi|)^n}{(A_0-A_{1}-A_{2})^n}+\epsilon, \forall t>T_{n,\epsilon}
\end{equation}

So we proceed by induction. Since $\epsilon>0,$ by Lemma \ref{LL1-1} and Lemma \ref{LL2}, There exists $T_{1,\epsilon}\ge 1$ { such that}
{
\begin{align*}
  1-\frac{C|\chi|}{A_0-A_{1}-A_{2}}-\frac{\epsilon }{M} \leq \underline{w}(t)&\leq w(t{+t_0},x)\cr
  &\le \overline{w}(t) \leq 1+\frac{C|\chi|}{A_0-A_{1}-A_{2}}+\frac{\epsilon }{M},\quad \forall t>T_{1,\epsilon},x\in \bar{\Omega}, {t_0\in\mathbb{R}}.
\end{align*}
}
Therefore 
\begin{align*}
     \frac{|u(t{+t_0},x)-u^*(t{+t_0},x)|}{|u^*(t{+t_0},x)|}=|w(t{+t_0},x)-1|\leq \frac{C|\chi|}{A_0-A_{1}-A_{2}}+\frac{\epsilon }{M},\quad  \forall t>T_{1,\epsilon},\ x\in \bar{\Omega},\ {t_0\in\mathbb{R}.}
\end{align*}
Thus 
{ 
\begin{align*}
|u(t{+t_0},x)-u^*(t{+t_0},x)|&\leq \frac{C|u^*(t{+t_0},x)||\chi|}{A_0-A_{1}-A_{2}}+\epsilon \cr
&\leq \frac{CM|\chi|}{A_0-A_{1}-A_{2}}+\epsilon,\quad \forall t>T_{1,\epsilon},\ x\in \bar{\Omega}, \ {t_0\in\mathbb{R}.} 
\end{align*}
}
So 
\begin{align}
\sup_{t_0 \in \mathbb{R}}\|u(t{+t_0},\cdot)-u^*(t{+t_0},\cdot)\|_\infty\leq \frac{CM|\chi|}{A_0-A_{1}-A_{2}}+\epsilon,\quad  \forall t>T_{1,\epsilon},
\end{align}
and the proof of equation \eqref{uniform-conv-eq00} follows for { $n = 1$}.

Next suppose that \eqref{uniform-conv-eq00} is true { up to some $n\ge 1$}, i.e {for all $0<\varepsilon<1$, there is $T_{n,\varepsilon}\ge n$ such that} 

$$\sup_{t_0 \in \mathbb{R}}\|u(t{+t_0},\cdot)-u^*(t{+t_0},\cdot)\|_\infty\leq M\frac{(C|\chi|)^n}{(A_0-A_{1}-A_{2})^n}+\epsilon, \forall t>T_{n,\epsilon},
$$
and we show that \eqref{uniform-conv-eq00} holds also for $n+1.$ Since equation \eqref{uniform-conv-eq00} holds for $n,$ { we get from} Lemma \ref{lem-02} that 
{ 
\begin{align*}
\label{gradients-bounds-eq1}
  &\quad {\sup_{t_0\in\mathbb{R}}}|\Delta(v^*-v)(t{+t_0})+\nabla \ln{u^*}(t{+t_0})\cdot\nabla(v^*-v)(t{+t_0})|\cr
  \leq& C_2\Big(M\frac{(C|\chi|)^n}{(A_0-A_{1}-A_{2})^n}+\epsilon+e^{-\tilde{\lambda}_0(t-t^*(u_0,v_0))}\Big)\cr
  \le& C_2\Big(M\frac{(C|\chi|)^n}{(A_0-A_{1}-A_{2})^n}+2\epsilon\Big):= \tilde{C}_{n,\varepsilon}\quad \forall \, t\ge T_{n,\epsilon}-\frac{\ln(\varepsilon)}{\tilde{\lambda}_0}+t^*(u_0,v_0).   
\end{align*}
}
Set $\tilde \epsilon={ \frac{\epsilon}{M}} (1-{ 2}\frac{ MC_2|\chi|}{A_0-A_{1}-A_{2}}).$ Note that $\tilde \epsilon>0$ since $MC_2=C$ and ${ 2}|\chi|<\frac{A_0-A_{1}-A_{2}}{C}$. Then, following the proof of Lemma \ref{LL1-1} and  Lemma \ref{LL2} with $C$ replaced by $\tilde{C}_{n,\varepsilon}$, there exists $T_{n+1,\tilde \epsilon}>n+1$ such that
\begin{equation*}
  1-\frac{\tilde{C}_{n,\varepsilon}|\chi|}{A_0-A_{1}-A_{2}}- \tilde \epsilon  \leq \underline{w}(t)\leq w(t{+t_0},x)\le \overline{w}(t) \leq 1+\frac{\tilde{C}_{n,\varepsilon}|\chi|}{A_0-A_{1}-A_{2}}+ \tilde \epsilon, \forall t>T_{n+1,{ \tilde \epsilon}},x\in \bar{\Omega},\  {t_0\in\mathbb{R}}.
\end{equation*}
Thus
\begin{align*}
\sup_{t_0 \in \mathbb{R}}\|u(t{+t_0},\cdot)-u^*(t{+t_0},\cdot)\|_{\infty} & \leq \frac{\tilde{C}_{n,\varepsilon} M|\chi|}{A_0-A_{1}-A_{2}}+\tilde{\epsilon} { \|u^*(t+t_0,\cdot)\|_{\infty} } \cr
&{ \leq} M\frac{(C|\chi|)^{n+1}}{(A_0-A_{1}-A_{2})^{n+1}}+{ 2}\frac{ MC_2|\chi|}{A_0-A_{1}-A_{2}}\epsilon+\tilde \epsilon{ M} \cr
&=M\frac{(C|\chi|)^{n+1}}{(A_0-A_{1}-A_{2})^{n+1}}+ \epsilon,\quad  \forall t>T_{n+1,\epsilon}.
\end{align*}
Therefore,  equation \eqref{uniform-conv-eq00} holds for all $n \in \mathbb{N}.$ It follows from Lemma \ref{lem-00} that  
\begin{align}
  \label{v-conv-eq0}
 \sup_{t_0 \in \mathbb{R}}\|v(t{+t_0},\cdot)-v^*(t{+t_0},\cdot)\|_\infty&\leq  
   C_0 \sup_{t_0 \in \mathbb{R}}\|u(t{+t_0},\cdot)-u^*(t{+t_0},\cdot)\|_\infty \cr
   & \leq 
   C_0\left(M\frac{(C|\chi|)^n}{(A_0-A_{1}-A_{2})^n}+\epsilon\right), \forall t>T_{n,\epsilon}.
\end{align}
Therefore the { desired result follows}, that is \eqref{global-stability-2-eq1}, from inequalities \eqref{uniform-conv-eq00} and equation \eqref{v-conv-eq0}, since $\frac{(C|\chi|)^n}{(A_0-A_{1}-A_{2})^n} \to 0 $ as $n \to \infty.$

{ Finally, we prove the uniqueness of $(u^*(t,x),v^*(t,x))$. To this end, suppose that $(\tilde{u}^*(t,x),\tilde{v}^*(t,x))$ is another positive entire solution of \eqref{u-v-eq00}.  Let $\varepsilon>0$. By \eqref{global-stability-2-eq1}, then there is $T_{\varepsilon}>0$ such that 
$$
\|\tilde{u}^*(t_0+t,\cdot)-u^*(t_0+t,\cdot)\|_{\infty}+\|\tilde{v}^*(t_0+t,\cdot)-v^*(t_0+t,\cdot)\|_{\infty}\le \varepsilon\quad \forall\ t\ge T_{\varepsilon},\ t_0\in\mathbb{R}.
$$
Hence, for any ${t}\in\mathbb{R}$, taking $t_0={t}-T_{\varepsilon}$, so that $t_0+T_{\varepsilon}=t$, we get 
\begin{align*}
 & \|\tilde{u}^*(t,\cdot)-u^*(t,\cdot)\|_{\infty}+\|\tilde{v}^*(t)-v^*(t,\cdot)\|_{\infty}\cr
 =&\|\tilde{u}^*(t_0+T_{\varepsilon},\cdot)-u^*(t_0+T_{\varepsilon},\cdot)\|_{\infty}+\|\tilde{v}^*(t_0+T_{\varepsilon},\cdot)-v^*(t_0+T_{\varepsilon},\cdot)\|_{\infty}\le \varepsilon.
\end{align*}
Since $\varepsilon>0$  is arbitrarily chosen, then $\|\tilde{u}^*(t,\cdot)-u^*(t,\cdot)\|_{\infty}+\|\tilde{v}^*(t)-v^*(t,\cdot)\|_{\infty}=0$  for all $t\in\mathbb{R}$, and hence $(u^*(t,\cdot),v^*(t,\cdot))=(\tilde{u}^*(t,\cdot),\tilde{v}^*(t,\cdot))$ for all $t\in\mathbb{R}$.
}

\end{proof}

{
\noindent {\bf Acknowledgment.} The authors also would like to thank the anonymous reviewers for their valuable comments and suggestions which considerably improved the presentation of this paper.
}

\begin{small}

\end{small}


\begin{thebibliography}{9}
\begin{small}



\bibitem{SAhmad1987}
 S. Ahmad, Convergence and ultimate bounds of solutions of the nonautonomous Volterra-Lotka competition equations, {\it J. Math. Anal. Appl.}, {\bf 127} (1987), no. 2, 377-387.

\vspace{-0.05 in}
 
\bibitem{ASV2009}
E. E. E. Arenas, A. Stevens, and  J. J. L. Velázquez, Simultaneous finite time blow-up in a two-species model for chemotaxis, {\it Analysis}, {\bf  29}(2009), 317-338 .


\vspace{-0.05 in}

\bibitem{NBYTMW05}
N. Bellomo, A. Bellouquid, Y. Tao, and M. Winkler, Toward a mathematical theory of Keller-Segel models of pattern formation in biological tissues, {\it Math. Models Methods Appl. Sci.}, {\bf 25(9)} (2015), 1663-1763.

\vspace{-0.05 in}

\bibitem{JDITBRB}
J. T. Doumatè, T. B. Issa, R. B. Salako, Competition-exclusion and coexistence in a two-strain SIS epidemic model in patchy environments.{ \it Discrete and Continuous Dynamical Systems} - B. (2023) doi: 10.3934/dcdsb.2023213

\vspace{-0.05 in}

\bibitem{Dan} D. Henry, Geometry theory of semilinear parabolic equations, Springer-Verlag, Berlin Heidelberg, New York 1981

\vspace{-0.05 in}

\bibitem{HVJ1997a}
M. A. Herrero and J. J. L. Velázquez,  A blow-up mechanism for a chemotaxis model,
{\it Ann. Sc. Norm. Super. Pisa, Cl. Sci.},  IV. Ser. 24 (1997), 633-683.

\vspace{-0.05 in}

\bibitem{HMVJ1996}
M. A. Herrero,  E.  Medina,  and J. J. L.  Velázquez,  Singularity patterns in a chemotaxis
model, {\it  Math. Ann.}, {\bf 306} (1996), 583-623.



\vspace{-0.10in}\bibitem{HeSh} G. Hetzer and W. Shen, Convergence in almost periodic competition diffusion systems,
{\it  J. Math. Anal. Appl.} {\bf 262} (2001),  307-338.


\vspace{-0.05 in}


\bibitem{THKJP09}
 T. Hillen, K. J. Painter, A users guide to PDE models for chemotaxis ,
{\it Math. Biol.}  {\bf  58} (2009) 183-217.

\vspace{-0.05 in}

\bibitem{H03}
D. Horstmann,  From 1970 until present: The Keller-Segel model in chemotaxis and its consequences, {\it  I. Jber. DMW},
{\bf  105} (2003), 103-165.

\vspace{-0.05 in}

\bibitem{HW2001}
D. Horstmann and G.  Wang,  Blow-up in a chemotaxis model without symmetry
assumptions, {\it  Eur. J. Appl. Math.}, {\bf  12}(2001), 159-177.

\vspace{-0.05 in}

\bibitem{ISM04}
M. Isenbach,  Chemotaxis. {\it Imperial College Press}, London (2004).

\vspace{-0.05 in}

\bibitem{ITBRS17}
T. B. Issa and  R. Salako, Asymptotic dynamics in a two-species chemotaxis model with non-local terms,
{\it  Discrete Contin. Dyn. Syst. Ser. B}, {\bf 22} (2017), no. 10, 3839-3874.

\vspace{-0.05 in}

\bibitem{ITBWS16}
T. B. Issa and  W. Shen, Dynamics in chemotaxis models of parabolic-elliptic type on bounded domain with time and space dependent logistic sources,
 {\it SIAM J. Appl. Dyn. Syst.}, {\bf 16} (2017), no. 2, 926-973.

\vspace{-0.05 in}

\bibitem{ITBWS17a}
T. B. Issa and  W. Shen, Persistence, coexistence and extinction in two species chemotaxis models on bounded heterogeneous environments,
{\it J. Dyn. Diff. Equat.}, (2018).  https://doi.org/10.1007/s10884-018-9686-7

\vspace{-0.05 in}

\bibitem{ITBWS17b}
T. B. Issa and  W. Shen, Uniqueness and stability of coexistence states in two species models with/without chemotaxis on bounded heterogeneous environments, {\it Journal of Dynamics and Differential Equations} (2018) https://doi.org/10.1007/s10884-018-9706-7.

\vspace{-0.05 in}

\bibitem{ITBWS2020}
T. B. Issa and  W. Shen, Pointwise persistence in  full chemotaxis models with logistic source on bounded heterogeneous environments, {\it Journal of Mathematical Analysis and Applications}, {\b 490(1)}(2020), 124204 .

\vspace{-0.05 in}

\bibitem{ITBRBSWS2021}
T. B. Issa, R.B. Salako, W. Shen, Traveling wave solutions for two species competitive chemotaxis systems, {\it Nonlinear Analysis}, { \bf 212} (2021).

\vspace{-0.05 in}

\bibitem{JaW92}
W. J$\ddot{a}$ger and S.  Luckhaus,   On explosions of solutions to a system of partial differential equations modeling chemotaxis, {\it Trans. Amer. Math. Soc.}, {\bf  329} (1992), 819-824.


\vspace{-0.05 in}


\bibitem{KS1970}
E. F. Keller and L. A.  Segel,   Initiation of slime mold aggregation viewed as an
instability, {\it  J. Theoret. Biol.}, {\bf  26} (1970),  399-415.

\vspace{-0.05 in}

\bibitem{KS71}
E. F. Keller and L. A. Segel,  A model for chemotaxis, {\it  J.Theoret. Biol.}, {\bf  30} (1971), 225-234.

\vspace{-0.05 in}

\bibitem{HIKWS2023b}
H. I. Kurt, W. Shen, Stabilization in two-species chemotaxis systems with singular sensitivity and Lotka-Volterra competitive kinetics. {\it Discrete and Continuous Dynamical Systems.} doi: 10.3934/dcds.2023130 

\vspace{-0.05 in}

\bibitem{HIKWS2023a}
H. I. Kurt, W. Shen, Chemotaxis systems with singular sensitivity and logistic source: Boundedness, persistence, absorbing set, and entire solutions, {\it Nonlinear Analysis: Real World Applications}, Volume 69, (2023),103762, ISSN 1468-1218.

\vspace{-0.05 in}

\bibitem{kuto_PHYSD}
   K. Kuto, K. Osaki,  T. Sakurai, and T. Tsujikawa,
   Spatial pattern formation in a chemotaxis-diffusion-growth model.
  {\it  Physica D}, {\bf 241} (2012), 1629-1639.

\vspace{-0.05 in}

\bibitem{DAL1991}
D. A. Lauffenburger, Quantitative studies of bacterial chemotaxis and microbial population dynamics
{\it Microbial. Ecol.}, {\bf 22}(1991),  175-85.

\vspace{-0.05 in}

\bibitem{LiKeWa} J. Li, Y. Ke, and Y. Wang,
 Large time behavior of solutions to a fully parabolic attraction-repulsion chemotaxis system with logistic source,
  {\it Nonlinear Anal. Real World Appl.} {\bf 39} (2018), 261-277.

  \vspace{-0.05 in}

  \bibitem{LiMu} K. Lin and C. Mu,
   Global dynamics in a fully parabolic chemotaxis system with logistic source,
   {\it  Discrete Contin. Dyn. Syst.} {\bf 36} (2016), no. 9, 5025-5046.

\vspace{-0.05 in}

\bibitem{NT1995}
T. Nagai,  Blow-up of radially symmetric solutions of a chemotaxis system,
{\it  Adv. Math.
Sci. Appl.}, {\bf  5} (198), 581-601.

\vspace{-0.05 in}

\bibitem{NT2001}
T. Nagai,  Blowup of nonradial solutions to parabolic-elliptic systems modeling
chemotaxis in two-dimensional domains,
{\it  J. Inequal. Appl.}, {\bf  6} (2001), 37-55.

\vspace{-0.05 in}

\bibitem{NT13}
M. Negreanu  and J. I. Tello, On a competitive system under chemotaxis effects with non-local terms,
{\it Nonlinearity},
{\bf 26} (2013), 1083-1103.


\vspace{-0.05 in}


\bibitem{NTa15}
M. Negreanu and J. I.  Tello, Asysmptotic stability of a two species chemotaxis system with non-diffusive chemoattractant,
{\it J. Differential Eq.}, {\bf  258} (2015),  1592-1617.

\vspace{-0.05 in}

\bibitem{PaHi} K. J.  Painter, and T. Hillen, Spatio-temporal chaos in a chemotaxis model.
 {\it  Physica D}, {\bf 240}(2011), 363-375


\vspace{-0.05 in}

\bibitem{RBSWS17a}
 R. B. Salako and W. Shen, Global existence and asymptotic behavior of classical solutions to a parabolic-elliptic chemotaxis system with logistic source on $R^N$,
{\it J. Differential Equations}, {\bf  262(11)} (2017), 5635-5690.

\vspace{-0.05 in}

\bibitem{RBSWS2018}
R. B. Salako, W Shen, Parabolic–elliptic chemotaxis model with space–time-dependent logistic sources on $\mathcal{R}^N.$ I. Persistence and asymptotic spreading,
{\it Mathematical Models and Methods in Applied Sciences} {\bf 28(11)} (2018), 2237-2273

\vspace{-0.05 in}

\bibitem{RBSWS2018b}
R. B. Salako, W. Shen, Parabolic–elliptic chemotaxis model with space–time dependent logistic sources on $\mathcal{R}^N.$ II. Existence, uniqueness, and stability of strictly positive entire solutions, {\it Journal of Mathematical Analysis and Applications}, {\bf 464 (1)} (2018), 883-910

\vspace{-0.05 in}

\bibitem{JMS1968}
J. M. Smith. “Mathematical Ideas in Biology,” Cambridge Univ. Press, London. 1968.

\vspace{-0.05 in}

\bibitem{T04}
J. I. Tello,  Mathematical analysis and stability of a chemotaxis problem with a logistic growth term,
{\it  Math. Methods  Appl. Sci.},  {\bf 27} (2004), 1865-1880.

\vspace{-0.05 in}


\bibitem{TJiCMAJ}
J. I. Tello and M. A. J.  Chaplain,   On the steady states and stability of solutions of a chemotaxis
problem with a logistic growth term,  Preprint.

\vspace{-0.05 in}

\bibitem{TW07}
J. I. Tello and W.  Winkler M, A chemotaxis system with logistic source, {\it Common Partial Diff. Eq.},
{\bf 32} (2007), 849-877.



\vspace{-0.05 in}


\bibitem{Win2010}
M. Winkler, Aggregation vs. global diffusive behavior in the higher-dimensional Keller-Segel model,
{\it  J.Differential Equations.}, {\bf 248} (2010),  2889-2905.

\vspace{-0.05 in}

\bibitem{W2011b}
M. Winkler,  Blow-up in a higher-dimensional chemotaxis system despite logistic growth restriction, {\it Journal of Mathematical Analysis and Applications}, {\bf 384}(2011), 261-272.

\vspace{-0.05 in}

\bibitem{Win2014}
M. Winkler, Global asymptotic stability of constant equilibria in a fully parabolic chemotaxis system with strong logistic damping,
{\it  J. Differential Equations.}, {\bf 257} (2014),  1056-1077.

\vspace{-0.05 in}

\bibitem{Zhe} J. Zheng,  Boundedness and global asymptotic stability of constant equilibria in a fully parabolic chemotaxis system with nonlinear logistic source,
    {\it  J. Math. Anal. Appl.} {\bf  450} (2017), no. 2, 1047-1061.

\end{small}

\end{thebibliography}
\end{document}